\documentclass[12pt,reqno]{amsart}
\usepackage{amsmath}
\usepackage{amsfonts}
\usepackage{amssymb}
\usepackage{amsthm}
\usepackage{amscd}
\usepackage{amsgen}
\usepackage{latexsym}
\usepackage{url}

\newtheorem{thm}{Theorem}[section]

\newtheorem{lemm}{Lemma}[section]

\newtheorem{rem}{Remark}[section]

\numberwithin{equation}{section}

\allowdisplaybreaks

\catcode`\@=11
\@namedef{subjclassname@2010}{%
  \textup{2010} Mathematics Subject Classification}
\catcode`\@=12

\def\po#1#2{(#1)_#2}

\newcommand{\Ric}{\operatorname{Ric}}
\newcommand{\scal}{\operatorname{scal}}
\newcommand{\tr}{\operatorname{tr}}

\def\S{{\mathbb S}}
\def\Rho{{\sf P}}               
\def\J{{\sf J}}                 
\def\f{{\frac{n}{2}}}
\def\M{{\mathcal M}}            
\def\st{\stackrel{\text{def}}{=}}

\title[Summation formulas for GJMS-operators and $Q$-curvatures]
{Summation formulas for GJMS-operators and $Q$-curvatures on the
M\"obius sphere}

\author{Andreas Juhl}

\address{Humboldt-Universit\"at, Institut f\"ur Mathematik, Unter den Linden,
D-10099 Berlin}
\email{ajuhl@math.hu-berlin.de}

\author{Christian Krattenthaler}

\address{Universit\"at Wien, Fakult\"at f\"{u}r Mathematik,
Nordbergstrasze 15, A-1090 Wien \newline \indent
\url{http://www.mat.univie.ac.at/~kratt/}}

\begin{document}

\dedicatory{Dedicated to Richard Askey}

\begin{abstract} For the M\"obius spheres $\S^{q,p}$, we give alternative
elementary proofs of the recursive formulas for GJMS-operators and
$Q$-curvatures due to the first author [{\it Geom.\ Funct.\ Anal.}
{\bf 23}, (2013), 1278--1370]. These
proofs make essential use of the theory of hypergeometric series.
\end{abstract}


\maketitle

\centerline \today

\tableofcontents

\renewcommand{\thefootnote}{}

\footnotetext{The work of the first author was supported by SFB 647
``Raum-Zeit-Materie'' of DFG. The work of the second author was
partially supported by the Austrian Science Foundation FWF, grants
Z130-N13 and S9607-N13, the latter in the framework of the National
Research Network ``Analytic Combinatorics and Probabilistic Number
Theory."}

\section{Introduction and statement of results}\label{intro}

In \cite{juhl-Q} and \cite{juhl-ex}, one of the authors proved
recursive formulas for GJMS-operators and Branson's $Q$-curvatures.
In the present article, we give alternative proofs of these results
for the space $\S^{q,p} = \S^q \times \S^p$ with the signature
$(q,p)$-metric $g_{\S^q} - g_{\S^p}$ given by the round metrics on
the factors.
The corresponding (explicit) formulas are found in
Theorems~\ref{op-form} and \ref{Q-form}. The proofs which we give
here provide another example of the extreme usefulness of
hypergeometric series identities. For other recent instances of
their crucial role in conformal differential geometry we refer to
\cite{GBR} and \cite{F-method}.

In the following, we write $n=q+p$.

A GJMS-operator $P_{2N}$, $N \ge 1$, is a specific rule which
associates to any pseudo-Riemannian manifold $(M,g)$ a differential
operator of the form
$$
P_{2N}(g) = \Delta_g^N + \mbox{lower-order terms},
$$
where $\Delta_g = -\delta_g d$ is the Laplace--Beltrami operator of
$g$. The operators $P_{2N}(g)$ are given by universal polynomial
formulas in terms of the metric $g$, its inverse, the curvature
tensor of $g$, and its covariant derivatives. Moreover, the operators
$P_{2N}(g)$ are conformally covariant in the sense that
\begin{equation}\label{covar}
e^{(\f+N)\varphi} P_{2N}(e^{2\varphi}g)(u) = P_{2N}(g)
\left(e^{(\f-N)\varphi}u \right)
\end{equation}
for all $\varphi \in C^\infty(M)$ and $u \in C^\infty(M)$. The
operators $P_{2N}(g)$ were derived in \cite{GJMS} from the $N^{th}$
powers of the Laplace--Beltrami operator of the Fefferman--Graham
ambient metric associated to $g$ \cite{FG-final}. On manifolds of
even dimension $n$, this construction is obstructed at $2N = n$,
i.e., it only yields a finite sequence $P_2(g), \dots, P_n(g)$ of
operators for which \eqref{covar} is valid for {\em arbitrary}
metrics. However, the obstructions should not be attributed to the
method of construction. In fact, Graham \cite{GJMS-2} proved that on
manifolds of dimension $4$, there is no conformally covariant cube
of the Laplacian. More generally, Gover and Hirachi \cite{GH} proved
that on manifolds of even dimension $n$, there is no conformally
covariant power of the Laplacian the order of which exceeds $n$.
Nevertheless, for even $n$ and locally conformally flat metrics, the
obstruction vanishes, and the construction of \cite{GJMS} still
yields an infinite sequence of conformally covariant operators. In
the present work, we will be concerned with such a case. For odd
$n$, the operators $P_{2N}(g)$ exist for general metrics $g$ and all
even orders $2N \ge 2$.

The GJMS-operators are not uniquely determined by the requirement of
their conformal covariance. It remains a challenge to characterize
them among other conformally covariant operators with leading part a
power of the Laplacian.

Explicit formulas for $P_2$ and $P_4$ are well-known. These
operators coincide with the Yamabe-operator
\begin{equation}\label{yam-op}
\Delta - \left(\f-1\right) \J, \quad \J = \frac{\scal}{2(n-1)},
\end{equation}
and the Paneitz-operator
\begin{equation}\label{pan-op}
\Delta^2 + \delta \left((n-2) \J g - 4 \Rho\right) d +
\left(\f-2\right) \left(\f \J^2 - 2 \Rho^2 - \Delta \J\right).
\end{equation}
Here
$$
\Rho = \frac{1}{n-2} \left(\Ric - \J g\right)
$$
denotes the Schouten tensor; in \eqref{pan-op}, it is regarded as an
endomorphism on one-forms.

Explicit formulas for higher-order GJMS-operators are much more
complicated. In \cite{juhl-book}, one of the authors introduced new
ideas for unveiling the structure of high-order GJMS-operators.
Among other things, these led to the discovery of a recursive
structure among these operators \cite{juhl-ex}.

Special cases of recursive formulas for GJMS-operators are the
relations
\begin{equation}\label{p4}
P_4 - P_2^2 = - 4 \delta(\Rho d) + \mu_4
\end{equation}
for the Paneitz operator $P_4$ and
\begin{equation}\label{p6}
P_6 - (2P_2 P_4 + 2P_4 P_2 - 3P_2^3) = - 48 \delta(\Rho^2 d) + \mu_6
\end{equation}
for the GJMS-operator $P_6$ of locally conformally flat metrics.
Here the scalar functions $\mu_4$ and $\mu_6$ are given in terms of
the curvature of $g$; for explicit formulas see \cite{juhl-ex}.

In order to describe the general recursive formula for
GJMS-operators, we need to introduce some notation. A sequence $I =
(I_1,\dots,I_r)$ of integers $I_j \ge 1$ will be regarded as a
composition of the sum $|I| = I_1 + I_2 + \cdots + I_r$. As usual, a
composition of a positive integer $m$ is a representation of $m$ as
a sum $m=m_1+m_2+\dots+m_r$ of positive integers
$m_1,m_2,\dots,m_r$, where two representations which contain the
same summands but differ in the order of the summands are regarded
as different. $|I|$ will be called the size of $I$. For any
composition $I$, we set
\begin{equation}
P_{2I} = P_{2I_1} \circ \cdots \circ P_{2I_r}.
\end{equation}
For $N \ge 1$, let
\begin{equation}\label{M-sum}
\M_{2N} \st \sum_{|I|=N} m_I P_{2I}
\end{equation}
with the multiplicities $m_I$ defined by
\begin{equation}\label{m-form}
m_I = -(-1)^r |I|! \, (|I|-1)! \prod_{j=1}^r \frac{1}{I_j!\,
(I_j-1)!} \prod_{j=1}^{r-1} \frac{1}{I_j + I_{j+1}}.
\end{equation}
Here, an empty product has to be interpreted as $1$. In particular,
$m_{(N)} = 1$ for all $N \ge 1$, and $m_I = m_{I^{-1}}$, where
$I^{-1}$ is the reverse composition of $I=(I_1,I_2,\dots,I_r)$, that
is, $I^{-1}=(I_r,I_{r-1},\dots,I_1)$.

Theorem 1.1 of \cite{juhl-ex} implies that, for locally conformally
flat Riemannian manifolds of dimension $n \ge 3$, we have
\begin{equation}\label{c1}
\M_{2N} = - c_N \delta (\Rho^{N-1} d) + \mbox{zeroth-order part},
\end{equation}
where
\begin{equation}\label{eq:cN}
c_N = 2^{N-1} N!\,(N-1)!.
\end{equation}

Eq.~\eqref{c1} may be regarded as a summation formula for the
operator $\M_{2N}$. We note that the sum \eqref{M-sum} contains
$2^{N-1}$ terms. Since one of these terms is $P_{2N}$, the operator
$\M_{2N}$ is of order $2N$. But \eqref{c1} shows that the operator
$\M_{2N}$ is actually only second-order. In other words, in the sum
\eqref{M-sum}, huge cancellations take place. Since $m_{(N)} = 1$,
Eq.~\eqref{c1} can be regarded as a formula for $P_{2N}^0 \st P_{2N}
- P_{2N}(1)$ in terms of a linear combination of the second-order
operator $\delta(\Rho^{N-1}d)$ and compositions of lower-order
GJMS-operators. In that sense, it is a recursive formula for
$P^0_{2N}$. Eqs.~\eqref{p4} and \eqref{p6} are the first two special
cases.

The proof of \eqref{c1} in \cite{juhl-ex} rests on the theory of
residue families as developed in \cite{juhl-book}. Residue families
are generalizations of certain differential intertwining operators
in representation theory (\cite{juhl-book}, \cite{KOSS}). It is the
recursive structure of these families which is responsible for the
identities \eqref{c1}. An alternative proof of \eqref{c1} was given
in \cite{FG-J}. It rests on the original definition of the
GJMS-operators in terms of the Fefferman--Graham ambient metric as
given in \cite{GJMS}. Note that the formulations in \cite{juhl-Q}
and \cite{juhl-ex} restrict to Riemannian metrics while in
\cite{FG-J} the signature of the metric is arbitrary.

Branson \cite{bran-2} used representation theory to derive explicit
expressions for the GJMS-operators of the spaces $\S^{q,p}$. In the
present paper, we use these formulas to give a direct proof of the
relation \eqref{c1} for $\S^{q,p}$. More precisely, we prove the
following summation formula.

\begin{thm}\label{op-form} On $\S^{q,p}$, we have
\begin{equation}\label{pseudo-even}
\M_{4N} = (2N)!\,(2N\!-\!1)!\, \left(\frac{1}{2}\!-\!B^2\!-\!C^2\right),
\; N \ge 1
\end{equation}
and
\begin{equation}\label{pseudo-odd}
\M_{4N+2} = (2N\!+\!1)!\,(2N)!\, (-B^2\!+\!C^2), \; N \ge 0.
\end{equation}
Here
\begin{equation}\label{BC-mobius}
B^2 = -\Delta_{\S^q} + \left(\frac{q-1}{2}\right)^2 \quad \mbox{and}
\quad C^2 = -\Delta_{\S^p} + \left(\frac{p-1}{2}\right)^2.
\end{equation}
\end{thm}

Theorem \ref{op-form} implies
\begin{align*}
\M_{4N}^0 & = (2N)!\,(2N\!-\!1)!\, (\Delta_{\S^q} + \Delta_{\S^p}), \; N \ge 1, \\
\M_{4N+2}^0 & = (2N\!+\!1)!\,(2N)!\, (\Delta_{\S^q} -
\Delta_{\S^p}), \; N \ge 0
\end{align*}
for the non-constant parts of the operators $\M_{2N}$. But using
\begin{equation}\label{rho}
\Rho = \frac{1}{2} \begin{pmatrix} 1_{\S^q} & 0 \\ 0 & -1_{\S^p}
\end{pmatrix},
\end{equation}
these two identities can be written in the form
\begin{align*}
\M_{4N}^0 & = - c_{2N} \delta (\Rho^{2N-1} d), \\
\M_{4N+2}^0 & = - c_{2N+1} \delta(\Rho^{2N} d),
\end{align*}
where $c_{2N}$ and $c_{2N+1}$ are defined in \eqref{eq:cN}. Thus,
Theorem~\ref{op-form} actually implies \eqref{c1}. Note that
Theorem~\ref{op-form} also determines the zeroth-order parts of all
$\M_{2N}$. It is easy to verify that the result fits with the
description of that part as given in \cite[Ex.~8.2]{juhl-ex}.

We continue with the discussion of $Q$-curvatures. The main result
of \cite{juhl-Q} concerns the zeroth-order parts of the
GJMS-operators. These define Branson's $Q$-curvatures as follows.
For even $n$, general metrics and $2N < n$, the equation
$$
P_{2N}(1) = (-1)^N \left( \f-N \right) Q_{2N}
$$
defines a scalar curvature quantity $Q_{2N}$. The quantities $Q_2,
\dots, Q_{n-2}$ are called the subcritical $Q$-curvatures. The
critical $Q$-curvature $Q_n$ can be defined through the subcritical
ones by continuation in dimension. In odd dimensions as well as for
locally conformally flat metrics in even dimensions, the
$Q$-curvatures $Q_{2N}$ are well defined for all $N\ge 1$. For
details we refer to \cite{bran-2} and \cite{juhl-book}.

Explicit formulas for $Q_2$ and $Q_4$ follow from \eqref{yam-op} and
\eqref{pan-op}:
$$
Q_2 = \J \quad \mbox{and} \quad Q_4 = \f \J^2 - 2 |\Rho|^2 - \Delta
\J.
$$
For $N \ge 3$, explicit formulas for $Q_{2N}$ are substantially more
complicated. Therefore, it is of some interest to establish
recursive formulas for $Q$-curvatures. Such a formula was found in
\cite{juhl-Q}.

In order to motivate its formulation, we start with the description
of two special cases. First of all, we rewrite $Q_4$ in the form
\begin{equation}\label{rec-q4}
Q_4 = - P_2(Q_2) - Q_2^2 + 2!\, 2^3\, v_4,
\end{equation}
where
$$
v_4 = \frac{1}{4} \tr (\wedge^2 \Rho).
$$
Next, for locally conformally flat metrics, we have (see
\cite{juhl-book})
\begin{equation}\label{rec-q6}
Q_6 = \left[-2P_2(Q_4) + 2 P_4(Q_2) - 3 P_2^2(Q_2)\right] - 6
\left[Q_4 + P_2(Q_2)\right] Q_2 - 2!\,3!\,2^5\, v_6,
\end{equation}
where
$$
v_6 = - \frac{1}{8} \tr (\wedge^3 \Rho).
$$
More conceptually, the coefficients $v_4$ and $v_6$ are Taylor
coefficients of the function
\begin{equation}\label{vol}
v(r) = \det \left(1-\frac{r^2}{4}\Rho\right) = 1 + v_2 r^2 + v_4 r^4
+ v_6 r^6 + \cdots,
\end{equation}
which describes the volume form of a Poincar\'e--Einstein metric
associated to the locally conformally flat metric $g$ (see
\cite{G-vol}, \cite{FG-final}). Eq.~\eqref{vol} implies
$$
v_{2k} = \left(-\frac{1}{2}\right)^k \tr (\wedge^k \Rho).
$$
Now, one can prove that the formulas \eqref{rec-q4} and
\eqref{rec-q6} are equivalent to the respective identities
\begin{equation}\label{Q4}
Q_4 + P_2(Q_2) = 2!\, 2^4\, w_4
\end{equation}
and
\begin{equation}\label{Q6}
-Q_6 - 2P_2(Q_4) + 2 P_4(Q_2) - 3 P_2^2(Q_2) = 2!\,3!\,2^6\, w_6,
\end{equation}
where the quantities $w_4$ and $w_6$ are Taylor coefficients of
$$
w(r) = \sqrt{v(r)},
$$
i.e.,
$$
\left(1+ r^2 w_2 + r^4 w_4 + r^6 w_6 + \cdots \right)^2 = v(r).
$$
More explicitly, we have
\begin{align*}
2 w_2 & = v_2, \\
2 w_4 & = \frac{1}{4} (4 v_4 - v_2^2), \\
2 w_6 & = \frac{1}{8} (8 v_6 - 4 v_4 v_2 + v_2^3).
\end{align*}

Now, in these terms, the main result of \cite{juhl-Q} implies that
\begin{equation}\label{c2}
\sum_{|I|=N} m_I \frac{P_{2I}(1)}{\f-I_{\text{last}}} = N!\,
(N\!-\!1)!\, w_{2N}
\end{equation}
for locally conformally flat metrics and $N \ge 1$. Here,
$I_{\text{last}}$ denotes the last entry of the composition $I$. For
$I=(J,\f)$, the quotient in \eqref{c2} is to be interpreted as
$$
(-1)^\f P_{2J}(Q_n).
$$
Again, the sum in \eqref{c2} contains $2^{N-1}$ terms. By separating
the contribution for $I=(N)$, Eq.~\eqref{c2} can be regarded as a
formula for $Q_{2N}$. It expresses $Q_{2N}$ in terms of lower-order
$Q$-curvatures, lower-order GJMS-operators and $w_{2N}$.
Eqs.~\eqref{Q4} and \eqref{Q6} are the first two special cases.

The proof of \eqref{c2} in \cite{juhl-Q} rests on the recursive
structure of residue families. This result played an important role
in the discussion of the recursive structure of GJMS-operators in
\cite{juhl-ex}. Again, an alternative proof of \eqref{c2} was given
in \cite{FG-J}.   

The following result proves \eqref{c2} for $\S^{q,p}$ by proving the
following summation formula.

\begin{thm}\label{Q-form} On $\S^{q,p}$, we have
\begin{equation}\label{Q-mobius}
\sum_{|I|=N} m_I \frac{P_{2I}(1)}{\f-I_{\text{last}}} = N!\,
(N\!-\!1)! \, \sum_{M=0}^N (-1)^M \binom{\frac{q}{2}}{M}
\binom{\frac{p}{2}}{N\!-\!M}
\end{equation}
for all $N \ge 1$.
\end{thm}

In order to see that Theorem \ref{Q-form} implies \eqref{c2}, we
note that \eqref{rho} and \eqref{vol} imply
$$
v(r) = (1-r^2/4)^q (1+r^2/4)^p.
$$
Hence
$$
w(r) = (1-r^2/4)^{q/2} (1+r^2/4)^{p/2}.
$$
It follows that the right-hand side of \eqref{Q-mobius} coincides
with
$$
N! \, (N-1)! \, 2^{2N}\, w_{2N}.
$$

The remaining part of the paper is organized as follows. In Section
\ref{op}, we prove Theorem \ref{op-form}. Section \ref{Q} contains
the proof of a summation formula which contains Theorem \ref{Q-form}
as a special case. The proofs of some technical results which are
used in the course of these proofs are collected in an appendix.

\section{Proof of Theorem \ref{op-form}}\label{op}

We start with a description of the GJMS-operators on $\S^{q,p}$.

\begin{thm}[\cite{bran-2}, Theorem 6.2]\label{b-mol} On $\S^{q,p}$,
the GJMS-operators are given by the product formulas
\begin{equation}\label{op-even}
P_{4N} = \prod_{j=1}^N
(B\!+\!C\!+\!(2j\!-\!1))(B\!-\!C\!-\!(2j\!-\!1))
(B\!+\!C\!-\!(2j\!-\!1))(B\!-\!C\!+\!(2j\!-\!1))
\end{equation}
and
\begin{equation}\label{op-odd}
P_{4N+2} = (-B^2\!+\!C^2) \prod_{j=1}^N
(B\!+\!C\!+\!2j)(B\!-\!C\!-\!2j)(B\!+\!C\!-\!2j)(B\!-\!C\!+\!2j).
\end{equation}
Here, the operators $B$ and $C$ are defined as the positive square
roots of the non-negative operators $B^2$ and $C^2$.
\end{thm}

The reader should note that the formulas \eqref{op-even} and
\eqref{op-odd} are equivalent to the product representations
\begin{equation}\label{op-even-b}
P_{4N} = \prod_{j=1}^{N}
\left((B^2\!-\!C^2)^2-2(2j\!-\!1)^2(B^2\!+\!C^2)+(2j\!-\!1)^4\right)
\end{equation}
and
\begin{equation}\label{op-odd-b}
P_{4N+2} = (-B^2\!+\!C^2) \prod_{j=1}^{N}
\left((B^2\!-\!C^2)^2 - 2(2j)^2(B^2\!+\!C^2)+(2j)^4 \right).
\end{equation}
An advantage of the latter formulas is that they do not require to
leave the framework of differential operators.

The formulas \eqref{op-even} and \eqref{op-odd} can be restated
uniformly as
\begin{equation}\label{poch1}
P_{2N} = 2^{2N} \left((C\!+\!B\!+\!1\!-\!N)/2 \right)_N
\left((C\!-\!B\!+\!1\!-\!N)/2 \right)_N,
\end{equation}
where $(\alpha)_m$ is the usual Pochhammer symbol defined by
$(\alpha)_m=\alpha(\alpha+1) \cdots(\alpha+m-1)$ for $m\ge1$, and
$(\alpha)_0=1$. In terms of the variables
\begin{equation}\label{XY}
X = C+B \quad \mbox{and} \quad  Y=C-B,
\end{equation}
the product formula \eqref{poch1} is equivalent to
\begin{equation}\label{eq:1}
P_{2N} = 2^{2N} \left((X+1-N)/2 \right)_N \, \left( (Y+1-N)/2
\right)_N.
\end{equation}
In the sequel, we shall regard $P_{2N}$ as this polynomial in the
variables $X$ and $Y$.

In the proofs of Theorems~\ref{op-form} and \ref{Q-form}, the basis
$\big(((X+1-A)/2)_A\big)_{A=0,1,\dots}$ of the linear space of all
polynomials in $X$ will play an essential role. (A glance at
\eqref{eq:1} may suggest why this could be the case.) In the following
lemma, we compute the structure coefficients with respect to
multiplication of this basis of polynomials. They will enter the
inductive proofs of Lemma~\ref{lem:2} and of Theorems~\ref{thm:1}
and \ref{thm:2}.

\begin{lemm}\label{lem:1} For all non-negative integers $A$ and $B$, we have
\begin{multline}\label{eq:2}
\left((X+1-A)/2\right)_A \left((X+1-B)/2 \right)_B \\
= \sum_{j=0}^{\lfloor(A+B)/2\rfloor} (-1)^j \frac{(-A/2)_j \,
(-B/2)_j \, (-(A+B)/2)_j}{j!} \left((X \!+ \!1\! - \!A\! - \!B\! +
\!2j)/2 \right)_{A+B-2j}.
\end{multline}
\end{lemm}

\begin{rem} The sum on the right-hand side of \eqref{eq:2}
terminates {\it at the latest} at $j={\lfloor(A+B)/2\rfloor}$. If
$A$ or $B$ should be even, it terminates already at $j=A/2$,
respectively at $j=B/2$.
\end{rem}

\begin{proof} The reader should recall the standard hypergeometric notation
$$
{}_p F_q\!\left[\begin{matrix} a_1,\dots,a_p\\
b_1,\dots,b_q\end{matrix}; z\right]=\sum _{m=0} ^{\infty}\frac
{\po{a_1}{m}\cdots\po{a_p}{m}} {m!\,\po{b_1}{m}\cdots\po{b_q}{m}}
z^m \,
$$
where $\po{a_1}{m}$, etc.\ are again Pochhammer symbols.
In terms of this notation, the right-hand side of \eqref{eq:2} reads
$$
\big((X+1-A-B)/2\big)_{A+B}\,
{}_{3} F_{2} \!\left [ \begin{matrix}
-(A+B)/2, -A/2, -B/2 \\
(1+X-A-B)/2,(1-X-A-B)/2 \end{matrix} ;
{\displaystyle 1}\right ].
$$
The $_3F_2$-series can be evaluated by means of the
Pfaff--Saalsch\"utz summation (cf.\ \cite[(2.3.1.3); Appendix
(III.2)]{SlatAC})
$$
{}_{3} F_{2} \! \left[ \begin{matrix} { a, b, -n} \\ { c, 1 + a + b
- c - n} \end{matrix} ; {\displaystyle 1} \right ] =
{\frac{({\textstyle c-a})_{n} \, ({\textstyle c-b})
_{n}}{({\textstyle c})_{n} \, ({\textstyle c-a-b})_{n}}},
$$
where $n$ is a non-negative integer. If we apply the formula, then
we obtain the left-hand side of \eqref{eq:2} after little
simplification.
\end{proof}

\begin{lemm}\label{lem:2} For all positive integers $a<N$, the
partial sum
\begin{equation}\label{part-sum}
S(N,a) = \sum_{J: \vert J\vert + a = N}^{} m_{(J,a)} P_{2J}
\end{equation}
satisfies
\begin{multline}\label{eq:3}
S(N,a) = \binom {N-1}{a-1} \sum_{k=0}^{\lfloor{(N-a)/2}\rfloor}
\sum_{l=0}^{\lfloor{(N-a)/2}\rfloor} (-1)^{N+k+l+a}2^{2N-2k-2l-2a} \\
\cdot \big((X+1-N+a+2k)/2\big)_{N-a-2k} \, \big((Y+1-N+a+2l)/2\big)_{N-a-2l}\\
\cdot \frac{(-N+a)_{2k}\,(-N+a)_{2l}\,(-N/2)_{k}\,(-N/2)_{l}} {k!\,l!}\,
{}_{4} F_{3} \!\left [\begin{matrix}
{-\frac{1}{2},-k,-l,\frac {1} {2}-\frac {N} {2}}\\
{-\frac{N}{2},\frac {a} {2}-\frac {N} {2},\frac {1} {2}+\frac {a}
{2}-\frac{N}{2}}
\end{matrix} ; {\displaystyle 1}\right ].
\end{multline}
\end{lemm}

\begin{proof} We prove the claim by induction on $N$, the case $N=2$ being
straightforward to check. Suppose that we have already established
\eqref{eq:3} up to $N-1$. If $a+\vert J\vert=N$, we write
$m_{(a,J)}$ as
\begin{equation*}
m_{(a,J)}=\begin{cases}
\displaystyle
-\frac {N!\,(N-1)!} {(N-a)!\,(N-a-1)!}
\frac {1} {N\cdot a!\,(a-1)!}m_{(N-a)},&\text{if }J=(N-a),\\
\displaystyle
-\frac {N!\,(N-1)!} {(N-a)!\,(N-a-1)!}
\frac {1} {(a+b)\,a!\,(a-1)!}m_{(b,K)},&\text{if }J=(b,K),
\end{cases}
\end{equation*}
and, if $J=(b,K)$, we write $P_{2J}=P_{2b}P_{2K}$.
Then the left-hand side of \eqref{eq:3} becomes
\begin{multline}\label{eq:3a}
- \frac {N!\,(N-1)!} {N\cdot a!\,(a-1)! \, (N-a)!\,(N-a-1)!}P_{2(N-a)}\\
- \sum_{b=1}^{N-a-1} \sum_{K:b+\vert K\vert=N-a}^{}
\frac{N!\,(N-1)!} {(a+b)\,a!\,(a-1)!\,(N-a)!\,(N-a-1)!}m_{(b,K)}
P_{2b} P_{2K}.
\end{multline}
If we now use the induction hypothesis, then this sum simplifies to
\begin{multline}\label{eq:3b}
-\sum_{b=1}^{N-a} \frac {N!\,(N-1)!}
{(a+b)\,a!\,(a-1)!\,(N-a)!\,(N-a-1)!} 2^{2b}\big((X+1-b)/2\big)_b\,\big((Y+1-b)/2\big)_b \\
\cdot \binom{N-a-1}{b-1}
\sum_{k=0}^{\lfloor{(N-a-b)/2}\rfloor}
\sum_{l=0}^{\lfloor{(N-a-b)/2}\rfloor}
(-1)^{N+k+l+a+b}2^{2N-2k-2l-2a-2b}\\
\cdot
\big((X+1-N+a+b+2k)/2\big)_{N-a-b-2k}\,
\big((Y+1-N+a+b+2l)/2\big)_{N-a-b-2l}\\
\cdot \frac {(-N+a+b)_{2k}\,(-N+a+b)_{2l}\,((-N+a)/2)_{k}\,((-N+a)/2)_{l}}
{k!\,l!} \\
\cdot {}_{4} F_{3} \! \left [ \begin{matrix}
{-\frac{1}{2},-k,-l,\frac{1}{2}+\frac{a}{2}-\frac{N}{2}}\\
{\frac{a}{2}-\frac{N}{2},\frac{a}{2}+\frac{b}{2}-\frac{N}{2},
\frac{1}{2}+\frac{a}{2}+\frac{b}{2}-\frac{N}{2}}
\end{matrix} ; {\displaystyle 1}\right ].
\end{multline}
Here, the term containing $P_{2(N-a)}$, appearing in \eqref{eq:3a}, arises as the summand
for $b=N-a$, since this forces $k$ and $l$ to be zero.

Now, Lemma \ref{lem:1} implies the expansions
\begin{multline*}
\big((X+1-b)/2\big)_b\,\big((X+1-N+a+b+2k)/2\big)_{N-a-b-2k}\\
= \sum_{{j_1}=0}^{\lfloor(N-a-2k)/2\rfloor}
(-1)^{j_1}\frac {(-b/2)_{j_1}\,(-(N-a-b-2k)/2)_{j_1}\,(-(N-a-2k)/2)_{j_1}}
{{j_1}!} \\
\cdot \big((X+1-N+a+2k+2{j_1})/2\big)_{N+a-2k-2{j_1}}.
\end{multline*}
and
\begin{multline*}
\big((Y+1-b)/2\big)_b\,\big((Y+1-N+a+b+2l)/2\big)_{N-a-b-2l}\\
= \sum_{{j_2}=0}^{\lfloor(N-a-2l)/2\rfloor}
(-1)^{j_2}\frac{(-b/2)_{j_2}\,(-(N-a-b-2l)/2)_{j_2}\,(-(N-a-2l)/2)_{j_2}}{{j_2}!} \\
\cdot \big((Y+1-N+a+2l+2{j_2})/2\big)_{N+a-2l-2{j_2}}.
\end{multline*}

We use these in \eqref{eq:3b} and, in addition, perform the index
transformation $s_1=k+j_1$ and $s_2=l+j_2$. Thus, we arrive at the
expression
\begin{multline}\label{eq:3c}
-\sum_{s_1=0}^{\lfloor(N-a)/2\rfloor}
\sum_{s_2=0}^{\lfloor(N-a)/2\rfloor}
\big((X+1-N+a+2s_1)/2\big)_{N-a-2s_1} \,
\big((Y+1-N+a+2s_2)/2\big)_{N-a-2s_2} \\
\cdot
\sum_{b=1}^{N-a}
\frac {N!\,(N-1)!} {(a+b)\,a!\,(a-1)!\,(N-a)!\,(N-a-1)!}
\binom {N-a-1}{b-1}
(-1)^{N+s_1+s_2+a+b}2^{2N-2a} \\
\cdot
\sum_{k=0}^{\lfloor{(N-a-b)/2}\rfloor}
\sum_{l=0}^{\lfloor{(N-a-b)/2}\rfloor}
(-b/2)_{s_1-k}\,(-b/2)_{s_2-l}\,
((-N+a)/2)_{s_1}\,
((-N+a)/2)_{s_2} \\
\cdot
((-N+a+b)/2)_{s_1}\,
((-N+a+b)/2)_{s_2}
\frac {((-N+a+b+1)/2)_{k}\,((-N+a+b+1)/2)_{l}}
{k!\,l!\,(s_1-k)!\,(s_2-l)!} \\
\cdot
{}_{4} F_{3} \!\left [ \begin{matrix}
{-\frac {1} {2},-k,-l,\frac {1} {2}+\frac {a} {2}-\frac {N} {2}}\\
{\frac {a} {2}-\frac {N} {2},\frac {a} {2}+\frac {b} {2}-\frac {N} {2},
\frac {1} {2}+\frac {a} {2}+\frac {b}
{2}-\frac {N} {2}}
\end{matrix} ; {\displaystyle 1}\right ].
\end{multline}
In this expression, we now concentrate on the terms involving the
summation index $k$ only:
\begin{equation}\label{eq:3d}
\sum_{k=0}^{s_1} \frac{(-b/2)_{s_1-k}\,(-(N-a-b-1)/2)_{k}\,(-k)_s} {k!\,(s_1-k)!},
\end{equation}
where $s$ stands for the summation index of the $_4F_3$-series in
\eqref{eq:6}. Because of the term $(-k)_s$ in the numerator of the
summand, we may start the summation at $k=s$ (instead of at $k=0$).
Hence, if we write this sum in hypergeometric notation, we obtain
$$
\frac {(-b/2)_{s_1-s}\,(-(N-a-b-1)/2)_{s}\,(-s)_s}{s!}  {}_{2} F_{1}
\!\left [ \begin{matrix}
-s_1+s, \frac{1}{2}+\frac{a}{2}-\frac{N}{2}+\frac{b}{2}+s \\
1+\frac{b}{2}-s_1+s
\end{matrix} ; {\displaystyle 1}\right ].
$$
This $_2F_1$-series can be evaluated by means of the Chu--Vandermonde
summation formula (see \cite[(1.7.7), Appendix (III.4)]{SlatAC}),
so that the sum in \eqref{eq:3d} equals
$$
\frac {(-s_1)_s\,(-(N-a-b-1)/2)_{s}\,((1+a-N)/2)_{s_1}}
{s_1! \,((1+a-N)/2)_{s}}.
$$
Clearly, an analogous computation can be done for the sum over $l$ in
\eqref{eq:3c}. Altogether, we see that \eqref{eq:3c} simplifies to
\begin{multline}\label{eq:3e}
-\sum_{s_1=0}^{\lfloor(N-a)/2\rfloor}
\sum_{s_2=0}^{\lfloor(N-a)/2\rfloor}
\big((X+1-N+a+2s_1)/2\big)_{N-a-2s_1}\,
\big((Y+1-N+a+2s_2)/2\big)_{N-a-2s_2}\\
\cdot
\sum_{b=1}^{N-a}
(-1)^{N+s_1+s_2+a+b}2^{2N-2s_1-2s_2-2a}
\frac{N!\,(N-1)!\,\binom {N-a-1}{b-1}}
{(a+b)\,a!\,(a-1)!\,(N-a)!\,(N-a-1)!}
\\
\cdot
\frac {(-N+a)_{2s_1}\,
(-N+a)_{2s_2}\,
((-N+a+b)/2)_{s_1}\,
((-N+a+b)/2)_{s_2}}
{s_1!\,s_2!}\\
\cdot
{}_{4} F_{3} \!\left [ \begin{matrix}
{-\frac{1}{2},-s_1,-s_2,\frac {1} {2}+\frac {a} {2}+
\frac{b}{2}-\frac {N} {2}}\\
{\frac{a}{2}-\frac {N} {2},\frac {1} {2}+\frac {a} {2}-\frac {N} {2},
\frac{a}{2}+\frac {b} {2}-\frac {N} {2}}
\end{matrix} ; {\displaystyle 1}\right ]\\
= -\sum_{s_1=0}^{\lfloor(N-a)/2\rfloor}
\sum_{s_2=0}^{\lfloor(N-a)/2\rfloor}
\big((X+1-N+a+2s_1)/2\big)_{N-a-2s_1}\,
\big((Y+1-N+a+2s_2)/2\big)_{N-a-2s_2}\\
\cdot
\sum_{b=1}^{N-a}
(-1)^{N+s_1+s_2+a+b}2^{2N-2s_1-2s_2-2a}
\frac {N!\,(N-1)!}
{(a+b)\,a!\,(a-1)!\,(N-a)!\,(b-1)!}
\\
\cdot
\sum_{s=0}^{s_1}
\frac{(-N+a)_{2s_1}\,
(-N+a)_{2s_2}\,
((-N+a+b+2s)/2)_{s_1-s}\,
((-N+a+b+2s)/2)_{s_2-s}}
{s_1!\,s_2!} \\
\cdot
\frac {(-1/2)_s\,(-s_1)_s\,(-s_2)_s\,(N-a-2s)!}
{s!\,(N-a)!\,(N-a-b-2s)!}.
\end{multline}
Next we concentrate on the terms involving the summation index $b$
only:
\begin{multline*}
\frac {1} {(N-a-2s-1)!}\sum_{b=1}^{N-a} \frac{(-1)^b} {a+b}
\binom {N-a-2s-1}{b-1}\\
\cdot
((-N+a+b+2s)/2)_{s_1-s}\,
((-N+a+b+2s)/2)_{s_2-s}.
\end{multline*}
By Lemma \ref{lem:diff2} with $X=-a-1$ and $M=N-a-2s-1$, this is
equal to
\begin{multline*}
\frac {(-1)^{N-a-2s-1}} {(-N+2s)_{N-a-2s}}((-N+2s)/2)_{s_1-s}
((-N+2s)/2)_{s_2-s}\\
+\chi(s_1=s_2=(N-a)/2)\cdot 2^{-N+a+2s}\\
=
-\frac {a!} {(N-2s)!}\frac {(-N/2)_{s_1}\,
(-N/2)_{s_2}}
{(-N/2)_{s}^2}+\chi(s_1=s_2=(N-a)/2)\cdot 2^{-N+a+2s},
\end{multline*}
where $\chi(\mathcal S)=1$ if $\mathcal S$ is true and
$\chi(\mathcal S)=0$ otherwise. If we substitute this in
\eqref{eq:3e}, then we obtain
\begin{multline*}
\sum_{s_1=0}^{\lfloor(N-a)/2\rfloor}
\sum_{s_2=0}^{\lfloor(N-a)/2\rfloor}
\big((X+1-N+a+2s_1)/2\big)_{N-a-2s_1}\,
\big((Y+1-N+a+2s_2)/2\big)_{N-a-2s_2}\\
\cdot (-1)^{N+s_1+s_2+a}2^{2N-2s_1-2s_2-2a}
\frac {(N-1)!} {(a-1)!\,(N-a)!} \\
\cdot \frac {(-N+a)_{2s_1}\, (-N+a)_{2s_2} \, (-N/2)_{s_1}\,(-N/2)_{s_2}}{s_1!\,s_2!} \\
\cdot \sum_{s=0}^{s_1} \frac {(-1/2)_s\,(-s_1)_s\,(-s_2)_s\,N! \,
(N-a-2s)!}{s! \, (N-2s)! \, (-N/2)_s^2 \, (N-a)!} \\
+ \chi(N-a\text { is even})\cdot \frac {N!\,(N-1)!}{a!\,(a-1)!\,((N-a)/2)!^2}
\kern4cm\\
\times \sum_{s=0}^{(N-a)/2} \frac {2^{-N+a+2s}(-1/2)_s \,
(-(N-a)/2)_s^2 \, (N-a-2s)!}{s!}.
\end{multline*}
Here, the first sum is, upon rewriting, exactly equal to the
right-hand side of \eqref{eq:3} (except that $s_1$ and $s_2$ took
over the role of $k$ and $l$). On the other hand, if we write the
second sum in hypergeometric notation, we obtain
\begin{multline*}
\chi(N-a\text { is even})\cdot
\frac {N!\,(N-1)!}
{a!\,(a-1)!\,((N-a)/2)!^2}\\
\times
\sum_{s=0}^{(N-a)/2}
\frac {2^{-N+a+2s}(-1/2)_s\,(-(N-a)/2)_s^2\,(N-a-2s)!}
{s!}\\
=\chi(N-a\text { is even})\cdot
\frac {N!\,(N-1)!}
{a!\,(a-1)!\,((N-a)/2)!^2}
{2^{-N+a}\,(N-a)!}
{}_{2} F_{1} \!\left [ \begin{matrix}
-\frac {1} {2},\frac {a} {2}-\frac {N} {2}\\
\frac {a} {2}-\frac {N} {2}+\frac {1} {2}
\end{matrix} ; {\displaystyle 1}\right ].
\end{multline*}
The $_2F_1$-series can be evaluated by the Chu--Vandermonde summation
formula, so that the above expression becomes
$$
\chi(N-a\text { is even})\cdot\frac {N!\,(N-1)!}
{a!\,(a-1)!\,((N-a)/2)!^2} {2^{-N+a}\,(N-a)!}  \frac
{((a-N+2)/2)_{(N-a)/2}} {((a-N+1)/2)_{(N-a)/2}},
$$
which vanishes because of the term $((a-N+2)/2)_{(N-a)/2}$. This
concludes the induction step. \end{proof}

\begin{thm}\label{thm:1} For all positive integers $N$, we have
\begin{equation}\label{eq:4}
\sum_{\vert I\vert=N}^{}m_IP_{2I} =
\begin{cases} N!\,(N-1)!\,XY,&\text {if $N$ is odd},\\
\frac {1}{2} N!\,(N-1)! \, \left(1-X^2-Y^2\right), & \text {if $N$
is even},
\end{cases}
\end{equation}
where $P_{2N}$ is defined in \eqref{eq:1}.
\end{thm}

\begin{proof} By the definition of $S(N,a)$ in \eqref{part-sum}, we have
\begin{equation}\label{eq:Ausdr}
\sum_{\vert I\vert=N}^{} m_I P_{2I} = P_{2N} + \sum_{a=1}^{N-1}
S(N,a) P_{2a}.
\end{equation}
If, in this equation, we substitute the expressions for $P_{2a}$ and
$S(N,a)$ given by \eqref{eq:1} and \eqref{eq:3}, then we see that
the left-hand side of \eqref{eq:4} is equal to
\begin{multline}\label{eq:5}
\sum_{a=1}^{N}
\big((X+1-a)/2\big)_{a}\,
\big((Y+1-a)/2\big)_{a}
\binom {N-1}{a-1}\\
\cdot
\sum_{k=0}^{\lfloor{(N-a)/2}\rfloor}
\sum_{l=0}^{\lfloor{(N-a)/2}\rfloor}
(-1)^{N+k+l+a}2^{2N-2k-2l}\kern7cm\\
\cdot
\big((X+1-N+a+2k)/2\big)_{N-a-2k}\,
\big((Y+1-N+a+2l)/2\big)_{N-a-2l}\\
\cdot
\frac {(-N+a)_{2k}\,(-N+a)_{2l}\,(-N/2)_{k}\,(-N/2)_{l}} {k!\,l!}\,
{}_{4} F_{3} \!\left [ \begin{matrix}
{-\frac {1} {2},-k,-l,\frac {1} {2}-\frac {N} {2}}\\
{-\frac {N} {2},\frac {a} {2}-\frac {N} {2},\frac {1} {2}+\frac {a}
{2}-\frac {N} {2}}
\end{matrix} ; {\displaystyle 1}\right ].
\end{multline}
Here, the term $P_{2N}$, appearing on the left-hand side of
\eqref{eq:4}, arises as the summand for $a=N$, since this forces $k$
and $l$ to be zero.

Lemma \ref{lem:1} implies the expansions
\begin{multline*}
\big((X+1-a)/2\big)_a\,\big((X+1-N+a+2k)/2\big)_{N-a-2k}\\
= \sum_{{j_1}=0}^{\lfloor(N-2k)/2\rfloor}
(-1)^{j_1}\frac {(-a/2)_{j_1}\,(-(N-a-2k)/2)_{j_1}\,(-(N-2k)/2)_{j_1}}
{{j_1}!} \\
\cdot
\big((X+1-N+2k+2{j_1})/2\big)_{N-2k-2{j_1}}.
\end{multline*}
and
\begin{multline*}
\big((Y+1-a)/2\big)_a\,\big((Y+1-N+a+2l)/2\big)_{N-a-2l}\\
= \sum_{{j_2}=0}^{\lfloor(N-2l)/2\rfloor}
(-1)^{j_2}\frac{(-a/2)_{j_2}\,(-(N-a-2l)/2)_{j_2}\,(-(N-2l)/2)_{j_2}}
{{j_2}!} \\
\cdot
\big((Y+1-N+2l+2{j_2})/2\big)_{N-2l-2{j_2}}.
\end{multline*}
We use these in \eqref{eq:5} and, in addition, perform the index
transformation $s_1=k+j_1$ and $s_2=l+j_2$. Thus, the
left-hand side in \eqref{eq:4} can be written in the form
\begin{multline}\label{eq:6}
\sum_{s_1=0}^{\lfloor N/2\rfloor} \sum_{s_2=0}^{\lfloor N/2\rfloor}
\sum_{a=1}^{N}(-1)^{N+s_1+s_2+a}2^{2N-2k-2l}
\binom {N-1}{a-1}\\
\cdot
\big((X+1-N+2s_1)/2\big)_{N-2s_1}
\big((Y+1-N+2s_2)/2\big)_{N-2s_2}
\\
\cdot
\sum_{k=0}^{\lfloor{(N-a)/2}\rfloor}
\sum_{l=0}^{\lfloor{(N-a)/2}\rfloor}
\frac{(-a/2)_{s_1-k}\,(-a/2)_{s_2-l}\,
(-(N-a)/2)_{s_1}\,
(-(N-a)/2)_{s_2}}
{(s_1-k)!\,(s_2-l)!\,
(-(N-a)/2)_{k}\,
(-(N-a)/2)_{l}
}\\
\cdot
\frac{(-N+a)_{2k}\,(-N+a)_{2l}\,(-N/2)_{s_1}\,(-N/2)_{s_2}} {k!\,l!}\,
{}_{4} F_{3} \!\left [ \begin{matrix}
{-\frac {1} {2},-k,-l,\frac {1} {2}-\frac {N} {2}}\\
{-\frac {N} {2},\frac {a} {2}-\frac {N} {2},\frac {1} {2}+\frac {a}
{2}-\frac {N} {2}}
\end{matrix} ; {\displaystyle 1}\right ].
\end{multline}
In this expression, we now concentrate on the terms involving the
summation index $k$ only:
\begin{multline}\label{eq:7}
\sum_{k=0}^{\lfloor{(N-a)/2}\rfloor}
2^{-2k}\frac {(-a/2)_{s_1-k}\,(-N+a)_{2k}\,(-k)_s}
{k!\,(s_1-k)!\,(-(N-a)/2)_{k}}\\=
\sum_{k=0}^{s_1}
\frac {(-a/2)_{s_1-k}\,(-(N-a-1)/2)_{k}\,(-k)_s}
{k!\,(s_1-k)!},
\end{multline}
where $s$ stands for the summation index of the $_4F_3$-series in
\eqref{eq:6}. Because of the term $(-k)_s$ in the numerator of the
summand, we may start the summation at $k=s$ (instead of at $k=0$).
Hence, if we write this sum in hypergeometric notation, we obtain
$$
\frac {(-a/2)_{s_1-s}\,(-(N-a-1)/2)_{s}\,(-s)_s} {s!}
{}_{2} F_{1} \!\left [ \begin{matrix}
-s_1+s,\frac {1} {2}-\frac {N} {2}+\frac {a} {2}+s \\
1+\frac {a} {2}-s_1+s
\end{matrix} ; {\displaystyle 1}\right ].
$$
This $_2F_1$-series can be evaluated by means of the Chu--Vandermonde
summation formula, so that the sum in \eqref{eq:7} equals
$$
\frac {(-s_1)_s\,(-(N-a-1)/2)_{s}\,((1-N)/2)_{s_1}}
{s_1! \,((1-N)/2)_{s}}.
$$
Clearly, an analogous computation can be done for the sum over $l$ in
\eqref{eq:6}. Altogether, we see that \eqref{eq:6} simplifies to
\begin{multline}\label{eq:8}
\kern-8pt
\sum_{s_1=0}^{\lfloor N/2\rfloor} \sum_{s_2=0}^{\lfloor N/2\rfloor}
\sum_{a=1}^{N}(-1)^{N+s_1+s_2+a}2^{2N}
\binom {N-1}{a-1}\\
\cdot
(-N/2)_{s_1}\,((1-N)/2)_{s_1}\,
(-N/2)_{s_2}\,((1-N)/2)_{s_2}
\\
\cdot
\big((X+1-N+2s_1)/2\big)_{N-2s_1}
\big((Y+1-N+2s_2)/2\big)_{N-2s_2}
\\
\cdot
\frac{
(-(N-a)/2)_{s_1}\,
(-(N-a)/2)_{s_2}}
{s_1!\,s_2!}
{}_{4} F_{3} \!\left [ \begin{matrix}
{-\frac {1} {2},-s_1,-s_2,\frac {1} {2}+\frac {a} {2}-\frac {N} {2}}\\
{-\frac {N} {2},\frac {1} {2}-\frac {N} {2},\frac {a}
{2}-\frac {N} {2}}
\end{matrix} ; {\displaystyle 1}\right ].
\end{multline}
Next we concentrate on the terms involving the summation index $a$
only:
\begin{align}
\notag
&\sum_{a=1}^{N}(-1)^a \binom {N-1}{a-1}
\frac {(-(N-a)/2)_{s_1}\,
(-(N-a)/2)_{s_2}\,(-(N-a-1)/2)_{s}}
{(-(N-a)/2)_{s}}\\
\notag
&\quad =
\frac {(N-1)!} {2^{2s}\,(N-2s-1)!}
\sum_{a=1}^{N-2s}(-1)^a \binom {N-2s-1}{a-1}\\
\notag
&\kern5cm\cdot
{(-(N-a-2s)/2)_{s_1-s}\,
(-(N-a-2s)/2)_{s_2-s}} \\
\notag
&\quad =
-\frac{(N-1)!} {2^{2s}\,(N-2s-1)!}
\sum_{a=0}^{N-2s-1}(-1)^a \binom {N-2s-1}{a}\\
&\kern5cm\cdot
{((a-N+2s+1)/2)_{s_1-s}\,
((a-N+2s+1)/2)_{s_2-s}},
\label{eq:8a}
\end{align}
where, in abuse of notation, we wrote again $s$ for the summation
index of the $_4F_3$-series in \eqref{eq:8}. By Lemma
\ref{lem:diff}, this sum vanishes if the degree of
$$
{((a-N+2s+1)/2)_{s_1-s}\, ((a-N+2s+1)/2)_{s_2-s}}
$$
as a polynomial in $a$ should be less than $N-2s-1$. In fact, this
is almost always the case since $s_1-s+s_2-s\le 2{\lfloor N/2\rfloor}-2s\le
N-2s$. More precisely, the sum in \eqref{eq:8a} is non-zero only if:
\begin{enumerate}
\item $N$ is odd and $s_1=s_2=(N-1)/2$;
\item $N$ is even and $s_1=s_2=N/2$ or
$\{s_1,s_2\}=\{N/2,(N-2)/2\}$.
\end{enumerate}
We now discuss these cases separately.
\smallskip

(1) Let $N$ be odd. Then, as we observed above, in the sum
\eqref{eq:8} only the terms corresponding to $s_1=s_2=(N-1)/2$
survive. Using Lemma \ref{lem:diff}, we see that \eqref{eq:8}
reduces to
\begin{multline*}
-4 \,(-N)_{N-1}^2 \, (X/2) \, (Y/2) \\
\times \sum_{s=0}^{(N-1)/2} (-1)^{N-2s}2^{-N+2s+1}\frac{(N-1)!}
{2^{2s}\,((N-1)/2)!^2} \frac {(-1/2)_s\,((1-N)/2)_s}
{s!\,(-N/2)_{s}} \\
= 2^{1-N} XY \frac{N!^2 \,(N-1)!}{((N-1)/2)!^2} {}_{2} F_{1} \!\left
[ \begin{matrix}
-\frac {1}{2},\frac{1}{2}-\frac{N}{2}\\
-\frac{N}{2}
\end{matrix} ; {\displaystyle 1}\right ].
\end{multline*}
The $_2F_1$-series can be evaluated by means of the Chu--Vandermonde
summation formula. After some simplification, we obtain
$N!\,(N-1)!\,XY$, in accordance with the claim in \eqref{eq:4}.
\smallskip

(2) Let $N$ be even. As we observed above, there are now three terms
in \eqref{eq:8} which survive: the terms for $s_1=s_2=N/2$, for
$(s_1,s_2)=(N/2,(N-2)/2)$, and for
$(s_1,s_2)=((N-2)/2,N/2)$.

We begin with the term for $s_1=N/2$ and $s_2=(N-2)/2$. By Lemma
\ref{lem:diff}, this term equals
\begin{multline*}
-4 \, (-N)_{N} \, (-N)_{N-2}\,((Y-1)/2)_2\\
\times \sum_{s=0}^{(N-2)/2} (-1)^{N-2s}2^{-N+2s+1}\frac{(N-1)!}
{2^{2s}\,(N/2)!\,((N-2)/2)!} \frac {(-1/2)_s\,((2-N)/2)_s}
{s!\,((1-N)/2)_{s}}\\
= -2^{-N}(Y^2-1) \frac{N!^2\,(N-1)!} {(N/2)!\,((N-2)/2)!} {}_{2}
F_{1} \! \left[ \begin{matrix} -\frac {1} {2},1-\frac {N} {2}\\
\frac{1}{2} - \frac {N}{2} \end{matrix} ; {\displaystyle 1} \right].
\end{multline*}
Again, the $_2F_1$-series can be evaluated by means of the
Chu--Vandermonde summation formula. After some simplification, we
obtain $-\frac {1} {2}N!\,(N-1)!\,(Y^2-1)$.  Clearly, an analogous
computation yields that the term for $s_1=(N-2)/2$ and $s_2=N/2$
equals $-\frac {1} {2}N!\,(N-1)!\,(X^2-1)$.  Finally, we consider
the term in \eqref{eq:8} corresponding to $s_1=s_2=N/2$. Here, in
order to apply Lemma \ref{lem:diff}, we need to compute the
coefficient $c_{N-2s-1}$ in the expansion
\begin{equation}\label{eq:exp}
{((a-N+2s+1)/2)_{(N-2s)/2}\, ((a-N+2s+1)/2)_{(N-2s)/2}} = \sum_{k=0}^{N-2s} c_k\binom ak.
\end{equation}
We have $c_{N-2s}=2^{-N+2s}(N-2s)!$. Comparison of coefficients of
$a^{N-2s-1}$ on both sides of \eqref{eq:exp} then yields that
\begin{align*}
c_{N-2s-1} & = 2^{-N+2s}\left(-2\sum_{\ell=1}^{(N-2s)/2}(2\ell-1)
+ \sum_{\ell=1}^{N-2s-1} \ell \right) \\
& = 2^{-N+2s}\left( -2 \left( \frac{N-2s}{2} \right)^2 +
\frac{(N-2s)(N-2s-1)}{2} \right) \\
& = - 2^{-N+2s-1}(N-2s).
\end{align*}
Hence, if we use this, together with Lemma \ref{lem:diff}, we obtain
that the term corresponding to $s_1=s_2=N/2$ in \eqref{eq:8} equals
\begin{multline*}
(-N)_{N}^2 \sum_{s=0}^{N/2}
(-1)^{N-2s+1}2^{-N+2s-1}(N-2s)\frac{(N-1)!} {2^{2s}\,(N/2)!^2}
\frac{(-1/2)_s\,(-N/2)_s}{s! \, ((1-N)/2)_{s}} \\
= -2^{-N-1} \frac{N!^2\,(N-1)!} {(N/2)!^2} \left( N \,{}_{2} F_{1}
\! \left[ \begin{matrix}
-\frac{1}{2},-\frac {N} {2}\\
\frac{1}{2}-\frac {N} {2}
\end{matrix} ; {\displaystyle 1}\right] +
\frac {N}{N-1}\,{}_{2} F_{1} \! \left [\begin{matrix}
\frac {1}{2},1-\frac{N}{2}\\
\frac {3}{2} - \frac{N}{2}
\end{matrix} ; {\displaystyle 1}\right ]\right).
\end{multline*}
Both $_2F_1$-series can be evaluated by means of the
Chu--Vandermonde summation formula. As it turns out, the first one
sums to zero.  If the result is simplified, one obtains $-\frac {1}
{2}N!\,(N-1)!$ eventually.

Putting our results together, we obtain
$$
\frac {1} {2}N!\,(N-1)!(-X^2+1-Y^2+1-1)=\frac {1} {2}N!\,(N-1)!
(1-X^2-Y^2)
$$
for the left-hand side of \eqref{eq:4}, which is in accordance with
the right-hand side. \end{proof}

Now, Theorem \ref{op-form} follows from Theorem \ref{thm:1} by
recalling \eqref{XY}.

\section{Proof of Theorem \ref{Q-form}}\label{Q}

As in Section \ref{op} we regard all GJMS-operators as polynomials
in the variables $X$ and $Y$ (see \eqref{eq:1}). In particular,
$P_{2N}$ is divisible by the monomial $(X+1-N)$.

\begin{thm}\label{thm:2} For all positive integers $N$, we have
\begin{multline}\label{eq:11}
\sum_{a=1}^{N} \sum_{\vert J\vert=N-a}^{}  m_{(J,a)} \frac{1}{X+1-a}
P_{2J} P_{2a} \\ = N!(N-1)! \sum_{k=0}^{N} (-1)^k
\binom{(X-Y+1)/2}{k} \binom{(X+Y+1)/2}{N-k}.
\end{multline}
\end{thm}

\begin{proof} By the definition of $S(N,a)$ in \eqref{part-sum}, we have
$$
\sum_{a=1}^{N} \sum_{\vert J\vert=N-a}^{}  m_{(J,a)} \frac{1}{X+1-a}
P_{2J} P_{2a} = \frac{1}{X+1-N} P_{2N} + \sum_{a=1}^{N-1}
\frac{1}{X+1-a} S(N,a) P_{2a}.
$$
The reader should note that the only difference to \eqref{eq:Ausdr}
are the factors $1/(X+1-N)$ in front of $P_{2N}$ and the factor
$1/(X+1-a)$ in front of the summand of the sum over $a$. We may
therefore proceed as in the proof of Theorem \ref{thm:1}, as long as
the summation over $a$ does not come into play. More precisely, by
comparing with \eqref{eq:8}, we see that the left-hand side of
\eqref{eq:11} simplifies to
\begin{multline}\label{eq:12}
\kern-8pt
\sum_{s_1=0}^{\lfloor N/2\rfloor} \sum_{s_2=0} ^{\lfloor N/2\rfloor}
\sum_{a=1}^{N} (-1)^{N+s_1+s_2+a}2^{2N}
\binom {N-1}{a-1}\\
\cdot
(-N/2)_{s_1}\,((1-N)/2)_{s_1}\,
(-N/2)_{s_2}\,((1-N)/2)_{s_2}
\\
\cdot \frac{1}{X+1-a} \big((X+1-N+2s_1)/2\big)_{N-2s_1}
\big((Y+1-N+2s_2)/2\big)_{N-2s_2} \\
\cdot \frac{ (-(N-a)/2)_{s_1}\, (-(N-a)/2)_{s_2}} {s_1!\,s_2!}
{}_{4} F_{3} \! \left[ \begin{matrix}
{-\frac{1}{2},-s_1,-s_2,\frac{1}{2} + \frac{a}{2} - \frac{N}{2}} \\
{-\frac{N}{2}, \frac{1}{2} - \frac{N}{2}, \frac{a}{2} - \frac{N}{2}}
\end{matrix} ; {\displaystyle 1}\right ].
\end{multline}
As in the proof of Theorem \ref{thm:1}, we now concentrate on the
terms involving the summation index $a$ only:
\begin{align*}
\notag & \sum_{a=1}^{N} (-1)^a \frac {1} {X+1-a}\binom {N-1}{a-1}
\frac {(-(N-a)/2)_{s_1}\, (-(N-a)/2)_{s_2}\,(-(N-a-1)/2)_{s}}
{(-(N-a)/2)_{s}}\\
\notag & \quad = \frac {(N-1)!} {2^{2s}\,(N-2s-1)!}
\sum_{a=1}^{N-2s}(-1)^a
\frac {1} {X+1-a} \binom {N-2s-1}{a-1}\\
\notag & \kern5cm\cdot {(-(N-a-2s)/2)_{s_1-s}\,
(-(N-a-2s)/2)_{s_2-s}} \\
\notag & \quad = -\frac {(N-1)!} {2^{2s}\,(N-2s-1)!}
\sum_{a=0}^{N-2s-1}(-1)^a
\frac {1} {X-a} \binom {N-2s-1}{a}\\
&\kern5cm \cdot {((a-N+2s+1)/2)_{s_1-s}\, ((a-N+2s+1)/2)_{s_2-s}},
\end{align*}
where, in abuse of notation, we wrote again $s$ for the summation
index of the $_4F_3$-series in \eqref{eq:12}. By Lemma
\ref{lem:diff2} with $M=N-2s-1$, this is equal to
\begin{multline}\label{eq:14}
\frac {(-1)^{N-2s}\,(N-1)!} {2^{2s}(X-N+2s+1)_{N-2s}}((X+1-N+2s)/2)_{s_1-s}
((X+1-N+2s)/2)_{s_2-s}\\
-\chi(s_1=s_2=N/2)\cdot 2^{-N}(N-1)!.
\end{multline}
We substitute this in \eqref{eq:12}. As before in the proof of Lemma
\ref{lem:2}, the second term in \eqref{eq:14} does not contribute
anything since the (remaining) sum over $s$ vanishes. Therefore,
substitution of \eqref{eq:14} in \eqref{eq:12} leads to the
expression
\begin{multline*}
\kern-8pt \sum_{s_1=0}^{\lfloor N/2\rfloor}\sum_{s_2=0}^{\lfloor
N/2\rfloor} \sum_{s=0}^{s_1} (-1)^{s_1+s_2+s}2^{2N-2s}
(-N/2)_{s_1}\,((1-N)/2)_{s_1} \,
(-N/2)_{s_2}\,((1-N)/2)_{s_2} \\
\cdot
\frac {\big((X+1-N+2s)/2\big)_{N-s-s_1}
\big((X+1-N+2s)/2\big)_{s_2-s}}
{(s_1-s)!\,(X-N+2s+1)_{N-2s}} \\
\cdot
\big((Y+1-N+2s_2)/2\big)_{N-2s_2}
\frac {(N-1)!\,(-1/2)_s\,(-s_2)_s}
{s_2!\,s!\,(-N/2)_s\,((1-N)/2)_s}.
\end{multline*}
We write the sum over $s_1$ in hypergeometric notation to obtain
\begin{multline*}
\sum_{s_2=0}^{\lfloor N/2\rfloor}
\sum_{s=0}^{s_1}
(-1)^{s_2}2^{2N-2s}
(-N/2)_{s_2}\,((1-N)/2)_{s_2}
\\
\cdot
\frac {\big((X+1-N+2s)/2\big)_{N-2s}
\big((X+1-N+2s)/2\big)_{s_2-s}}
{(X-N+2s+1)_{N-2s}}
\\
\cdot
\big((Y+1-N+2s_2)/2\big)_{N-2s_2}
\frac {(N-1)!\,(-1/2)_s\,(-s_2)_s}
{s_2!\,s!}
{}_{2} F_{1} \!\left [ \begin{matrix}
-\frac {N} {2}+s,\frac {1} {2}-\frac {N} {2}+s \\
-\frac {X} {2}+\frac {1} {2}-\frac {N} {2}+s
\end{matrix} ; {\displaystyle 1}\right ].
\end{multline*}
The $_2F_1$-series can be evaluated by means of the Chu--Vandermonde
summation formula. After some simplification, we arrive at
\begin{multline*}
\sum_{s_2=0}^{\lfloor N/2\rfloor}
\sum_{s=0}^{s_1}
(-1)^{s_2}2^{N}
(-N/2)_{s_2}\,((1-N)/2)_{s_2}\\
\cdot
\big((X+1-N+2s)/2\big)_{s_2-s}
\big((Y+1-N+2s_2)/2\big)_{N-2s_2}
\frac {(N-1)!\,(-1/2)_s\,(-s_2)_s}
{s_2!\,s!}\\
=\sum_{s_2=0}^{\lfloor N/2\rfloor}
(-1)^{s_2}2^{N}\,\frac {(N-1)!} {s_2!}
(-N/2)_{s_2}\,((1-N)/2)_{s_2}\kern3cm\\
\cdot
\big((X+1-N)/2\big)_{s_2}
\big((Y+1-N+2s_2)/2\big)_{N-2s_2}
\,{}_{2} F_{1} \!\left [ \begin{matrix}
-\frac {1} {2},-s_2\\
\frac {X} {2}+\frac {1} {2}-\frac {N} {2}
\end{matrix} ; {\displaystyle 1}\right ].
\end{multline*}
After having evaluated the $_2F_1$-series by means of the
Chu--Vandermonde summation formula, we are left with
\begin{multline*}
\sum_{s_2=0}^{\lfloor N/2\rfloor}
(-1)^{s_2}2^{N}\,\frac {(N-1)!} {s_2!}
(-N/2)_{s_2}\,((1-N)/2)_{s_2}\\
\cdot
\big((X+2-N)/2\big)_{s_2}
\big((Y+1-N+2s_2)/2\big)_{N-2s_2}\\
=2^{N}{(N-1)!}
\big((Y+1-N)/2\big)_{N}
\, {}_{3} F_{2} \! \left [ \begin{matrix}
\frac {X} {2}-\frac {N} {2}+1,-\frac {N} {2},\frac {1} {2}-\frac {N}{2} \\
-\frac {Y} {2}-\frac {N} {2}+\frac {1} {2},
\frac {Y} {2}-\frac {N} {2}+\frac {1} {2}
\end{matrix} ; {\displaystyle 1}\right ].
\end{multline*}
By the transformation formula in Lemma \ref{lem:3}, this is equal to
$$
(N-1)!\,((X+Y-2N+3)/2)_N \,
{}_{2} F_{1} \!\left [ \begin{matrix} - \frac{X}{2} + \frac{Y}{2} - \frac{1}{2},-N\\
\frac{X}{2}+\frac{Y}{2}-N+\frac{3}{2}
\end{matrix} ; {\displaystyle -1}\right ].
$$
which agrees with the right-hand side of \eqref{eq:11}. This
finishes the proof.
\end{proof}

Now we are finally in the position to establish Theorem~\ref{Q-form}.

\begin{proof}[Proof of Theorem \ref{Q-form}] Using $n=q+p$, we find
$$
X(1)= (C+B)(1) = \f-1 \quad \mbox{and} \quad Y(1)=(C-B)(1) =
\frac{p-q}{2}.
$$
Thus, we can write the constant term of $P_{2N}$ in the form
$$
P_{2N}\left(\f-1,\frac{p-q}{2}\right).
$$
Now we calculate
\begin{align*}
\sum_{|I|=N} m_I & \frac{1}{\f-I_{\text{last}}} P_{2I}(1) \\
& = \sum_{a=1}^N \sum_{|J|=N-a} m_{(J,a)} \frac{1}{\f-a}
P_{2J}\left(\f-1,\frac{p-q}{2}\right)
P_{2a}\left(\f-1,\frac{p-q}{2}\right) \\
& = N!(N-1)! \sum_{k=0}^N (-1)^k \binom{\frac{q}{2}}{k}
\binom{\frac{p}{2}}{N-k}.
\end{align*}
In the last step we have applied Theorem \ref{thm:2}.
\end{proof}


\appendix
\label{app}
\global\def\theTheorem{\mbox{A
}}
\setcounter{thm}{0}

\def\thesection{\!\!\!}
\section{Proofs of some auxiliary results}

\def\thesection{A}

Here we prove three technical lemmas which were used in Sections
\ref{op} and \ref{Q}.

\begin{lemm}\label{lem:diff} Let $p(x)$ be a polynomial in $x$, and
suppose that
$$
p(x) = \sum_{k\ge0}^{} c_k \binom{x}{k},
$$
for {\em(}uniquely determined\/{\em)} coefficients $c_k$.
Furthermore, let $M$ be an integer.  Then
\begin{equation*}
\sum_{a=0}^{M} (-1)^a \binom{M}{a} p(a) = (-1)^M c_M.
\end{equation*}
\end{lemm}

\begin{proof} This is a classical fact from finite difference
calculus. For the convenience of the reader, we give the simple
proof. We calculate
\begin{align*}
\sum_{a=0}^{M}(-1)^a \binom{M}{a} p(a)&=
\sum_{a=0}^{M}(-1)^a \binom{M}{a} \sum_{k\ge0}^{} c_k \binom{a}{k} \\
&=\sum_{k\ge0}^{} c_k \binom{M}{k}
\sum_{a=0}^{M} (-1)^a \binom{M-k}{a-k}\\
&=\sum_{k\ge0}^{} c_k \binom{M}{k} (-1)^k \delta_{M,k} = (-1)^M c_M.
\end{align*}
The proof is complete. \end{proof}

\begin{lemm}\label{lem:diff2} Let $p(x)$ be a polynomial in $x$, and
suppose that
$$
p(x)=\sum_{k\ge0}^{} c_k \binom{x}{k}
$$
for {\em(}uniquely determined\/{\em)} $c_k$. Furthermore, let $M$ be
an integer. If the degree of $p(x)$ is at most $M$, then
\begin{equation*}
\sum_{b=0}^{M} \frac{(-1)^b} {X-b} \binom{M}{b} p(b) = (-1)^M
\frac{M!}{(X-M)_{M+1}} p(X).
\end{equation*}
If the degree of $p(x)$ equals $M+1$, then
\begin{equation*}
\sum_{b=0}^{M} \frac {(-1)^b} {X-b} \binom {M}{b} p(b) = (-1)^M
\frac {M!}{(X-M)_{M+1}} p(X) + \frac {c_{M+1}}{M+1}.
\end{equation*}
\end{lemm}

\begin{proof} It suffices to verify the claim for a basis of the
vector space of polynomials in $x$ of degree at most $M+1$. We
choose $\left\{\binom{x}{k}\right\}_{0 \le k \le M+1}$ as such a
basis. If $k\le M$, we have
\begin{align*}
\sum_{b=0}^{M}\frac {(-1)^b} {X-b} \binom {M}{b} \binom{b}{k} & =
\frac{(-1)^k} {X-k} \binom{M}{k} {}_{2} F_{1} \!\left
[\begin{matrix} -M+k,k-X\\k-X+1 \end{matrix}
; {\displaystyle 1} \right ]\\
& = \frac{(-1)^k}{X-k}
\binom{M}{k} \frac{(M-k)!}{(-X+k+1)_{M-k}}\\
& =(-1)^M \frac {M!}{k!} \frac{(X-k+1)_k}{(X-M)_{M+1}}\\
& =(-1)^M \frac {M!}{(X-M)_{M+1}}\binom {X}{k},
\end{align*}
while if $k=M+1$ we obtain $0$. Here, we used the Chu--Vandermonde
summation formula to evaluate the $_2F_1$-series. All this is in
agreement with our claims.
\end{proof}

\begin{lemm}\label{lem:3} For all non-negative integers $N$, we have
\begin{equation}\label{eq:16}
{}_{3} F_{2} \! \left[ \begin{matrix} {a},
{\frac{1}{2}-{\frac{N}{2}}}, -{\frac{N}{2}} \\ {e,1-N-e}
\end{matrix} ; {\displaystyle 1} \right ] =  2^{-N}
\frac{(e-a)_N}{(e)_N} {}_{2}
F_{1} \!\left [ \begin{matrix} {1-a-e-N,-N} \\
{1+a-e-N} \end{matrix} ; {\displaystyle -1}\right ].
\end{equation}
\end{lemm}

\begin{proof} In \cite[(3.14)]{KrSrAB},
\begin{equation*}
{}_{4} F_{3} \!\left [ \begin{matrix} {a,b,{\frac 1 2} - {\frac
N 2}, -{\frac N  2}} \\
{{\frac 1 2} + a + b, e,1 -N - e}\end{matrix} ; {\displaystyle  1}\right ] =
\frac {(a+e)_N} {(e)_N}
{}_{3} F_{2} \!\left [ \begin{matrix} { 2\,a, a + b, -N}\\ { 2\,a + 2\,b, a + e} \end{matrix}
; {\displaystyle 1}\right ],
\end{equation*}
we let $b$ tend to infinity. As a result, we obtain the transformation
formula
\begin{equation*}
{}_{3} F_{2} \!\left [ \begin{matrix} { a, {\frac 1 2} - {\frac N 2},
-{\frac N 2}} \\
{ e,1 -N - e}\end{matrix} ; {\displaystyle  1}\right ]  =
\frac {(a+e)_N} {(e)_N}
{}_{2} F_{1} \!\left [ \begin{matrix} { 2\,a, -N}\\
{a + e}\end{matrix} ; {\displaystyle \frac {1} {2}}\right ].
\end{equation*}
If we now apply the transformation formula (see \cite[(1.7.1.3)]{SlatAC})
$$
{}_{2} F_{1}
\!\left [ \begin{matrix} { A,B}\\ { C}\end{matrix} ; {\displaystyle  z}\right ]
= {{{\left( 1 - z \right) }^{-B}}}
{{{}_{2} F_{1} \!\left [ \begin{matrix} {C-A,B} \\
{C} \end{matrix} ; {\displaystyle -{\frac z {1-z}}}\right ] }}
$$
to the $_2F_1$-series on the right-hand side, then we obtain
\begin{equation*}
2^{-N} \frac {(a+e)_N} {(e)_N}
{}_{2} F_{1} \!\left [ \begin{matrix} {e-a,-N} \\
{a + e}\end{matrix} ; {\displaystyle -1}\right ].
\end{equation*}
Finally, if we reverse the order of summands in this $_2F_1$-series
(that is, if we denote the summation index in the $_2F_1$-series by
$s$, then we replace $s$ by $N-s$ and rewrite the result again in
hypergeometric notation), then we arrive at the right-hand side of
\eqref{eq:16}. \end{proof}


\end{document}